\documentclass[12pt]{amsart}
\newtheorem{theorem}{Theorem}[section]
\newtheorem{lemma}[theorem]{Lemma}
\newtheorem{corollary}[theorem]{Corollary}
\theoremstyle{definition}

\newtheorem{example}[theorem]{Example}
\theoremstyle{remark}
\newtheorem{remark}[theorem]{Remark}
\numberwithin{equation}{section}

\newsymbol\subsetneq 2328
\newsymbol\dotplus 1275

\def\A{\bar A}
\def\z*{\bar z}

\def\B{\mathsf B}
\def\L{\mathsf L}
\def\uno{\mathsf 1}
\def\H{\mathcal H}
\def\fh{\mathfrak h}
\def\HM{\H_{-}}
\def\HP{\H_{+}}
\def\K{\mathcal K}
\def\X{\mathcal X}
\def\Y{\mathcal Y}
\def\N{\mathcal N}
\def\MG{\mathcal G}
\def\RE{\mathbb R}
\def\C{{\mathbb C}}
\def\NA{\mathbb N}
\def\G*{G_\star}
\def\ph*{\phi_\star}
\def\vp{\varphi}
\def\wtilde{\widetilde}

\begin{document}

\title[Self-Adjoint Extensions by Additive Perturbations]{Self-Adjoint Extensions by Additive Perturbations}

\author{Andrea Posilicano}

\address{Dipartimento di Scienze, Universit\`a dell'Insubria, I-22100
Como, Italy}

\email{andreap@uninsubria.it}

\begin{abstract}
Let $A_\N$ be the symmetric operator given by the restriction of $A$
to $\N$, where $A$ is a
self-adjoint operator on the Hilbert space $\H$ and $\N$ is a linear dense set which is 
closed with respect to the graph
norm on $D(A)$, the operator domain of $A$. 
We show that any self-adjoint extension $A_\Theta$ of
$A_\N$ such that $D(A_\Theta)\cap D(A)=\N$ can be additively
decomposed by the sum 
$A_\Theta=\A+T_\Theta$, 
where both the
operators $\A$ and $T_\Theta$ take values in the strong dual of
$D(A)$. The operator $\A$ is the closed extension of $A$ to the whole
$\H$ whereas $T_\Theta$ 
is explicitly written in terms of a (abstract) boundary condition depending on
$\N$ and on the extension parameter $\Theta$, a self-adjoint operator on 
an auxiliary Hilbert space isomorphic (as a set) to the
deficiency spaces of $A_\N$. The explicit connection with both Kre\u\i n's
resolvent formula and von
Neumann's theory of self-adjoint extensions is given. 
\end{abstract}

\maketitle

\section{Introduction}
Given a self-adjoint operator $A:D(A)\subseteq\H\to\H$, let $A_\N$ be
the restriction of $A$ to $\N$, where $\N\subsetneq D(A)$ is a dense linear subspace
which is closed with respect to the graph norm. Then $A_\N$
is a closed, densely defined, symmetric operator. Since $\N\neq D(A)$, 
$A_\N$ is not essentially
self-adjoint, as $A$ is a non-trivial extension of $A_\N$, and, by the
famed von Neumann's formulae \cite{[N]}, we know that
$A_\N$ has an infinite family of self-adjoint
extensions $A_U$ parametrized by the unitary maps $U$ from $\K_+$ onto $\K_-$,
where $\K_\pm:=$Kernel$\,(-A_\N^*\pm i)$ denotes the deficiency
spaces.\par 
In section 2 we define a family $A_\Theta$ of extensions of $A_\N$ by
means of a Kre\u\i n-like
formula i.e. by explicitly giving its resolvent
$(-A_\Theta+z)^{-1}$ (see Theorem 2.1). By using the approach developed in
\cite{[P1]}, we describe the domain of $A_\Theta$ in terms of
the boundary condition $\tau\ph*=\Theta\, Q_\phi$, where $\tau:
D(A)\to\fh$ is a surjective continuos linear mapping with Kernel$\,\tau=\N$,
$\Theta:D(\Theta)\subseteq\fh\to\fh$ is self-adjoint and $\fh$ is a
Hilbert space isomorphic (as a set) to $\K_\pm$. \par 
In section 3 we use the resolvent $(-A_\Theta+z)^{-1}$ given in
Theorem 2.1 to re-write $A_\Theta$ in a more appealing way as a sum $\A+T_\Theta$ 
where both $\A$ and $T_\Theta$ take values in the strong dual
(with respect to the graph norm) of $D(A)$ (see Theorem 3.1); $\A$ is nothing
else that the closed extension of $A$  to the whole
Hilbert space $\H$ and $T_\Theta$ is explicitly given in terms of the maps
$\tau$ and $\Theta$ giving the boundary conditions. This result gives
an extension, and a rephrasing in
terms of boundary conditions, of 
the results obtained in \cite{[K]} (and references therein, in
particular \cite{[KY]}), where 
$A$ is strictly positive and $\N$ is closed in $D(A^{1/2})$ 
(see Remark 3.5). As regards
boundary conditions the reader is also refered to \cite{[Ko]}, 
where $A=-\Delta+\lambda$,
$\lambda>0$, $\N$ the kernel of the evaluation map along a
regular submanifold, and to \cite{[P2]}, where $A$ is an arbitrary 
injective self-adjoint operator. 
\par Successively, is section 4, we study the connection of the
self-adjoint extensions defined in the previuos sections with the ones 
given by von Neumann's theory \cite{[N]}. We prove (see Theorem 4.1) that the
operator $\wtilde A=\A+T$
defined in Theorem 3.4, of which the self-adjoint
$A_\Theta=\A+T_\Theta$ is a restriction, coincides 
with $A^*_\N$; moreover we explicitly define a map on self-adjoint 
operators $\Theta:D(\Theta)\subseteq\fh\to\fh$ to unitary operators 
$U:\K_+\to
\K_-$ such that $A_\Theta=A_U$, where $A_U$ denotes the von
Neumann's  extension corresponding to $U$. Such correspondence is then
explicitly inverted (see Theorem 4.3). This shows (see Corollary 4.4) that
$\wtilde A=\A+T$ 
coincides with a self-adjoint extension $\widehat A$ of $A_\N$ such that
$D(\widehat A)\cap D(A)=\N$ if and only if the boundary
condition $\tau\ph*=\Theta\,Q_\phi$ holds for some self-adjoint operator
$\Theta$. 
\par In section 5 we conclude
with some examples both in the case of finite and infinite deficiency
indices. Example 5.1 (also see Remark 4.2) shows that, in the case 
dim$\,\K_\pm<+\infty$, our results reproduce the
theory of finite rank
perturbations as given in \cite{[AK]}, \S 3.1, and thus
they can be viewed as an extension of such a theory to the infinite rank
case. In example 5.2 we give two examples in the infinite rank case:
infinitely many point interaction in three dimensions and singular
perturbations, supported on $d$-sets with $0<n-d<2s$, of traslation
invariant pseudo-differential operators with domain the Sobolev space 
$H^s(\RE^n)$.  
\section*{Notations and definitions}
\begin{itemize}
\item Given a Banach space $\X$ we denote by $\X'$ its strong dual.  
\item $\L(\X,\Y)$ denotes the
space of linear operators 
from the Banach space $\X$ to the Banach space $\Y$; 
$\L(\X)\equiv\L(\X,\X)$. 
\item $\B(\X,\Y)$ denotes the Banach space of 
bounded, everywhere defined, linear operators 
on the Banach space $\X$ to the Banach space $\Y$; $\B(\X)\equiv\B(\X,\X)$. 
\item Given $A\in \L(\X,\Y)$ densely defined, the
closed operator $A'\in \L(\Y',\X')$ is the adjoint of $A$ 
i.e. 
$$\forall\,\phi\in D(A)\subseteq\X,\quad\forall
\lambda\in D(A')\subseteq
\Y',\qquad (A'\lambda)(\phi)=\lambda(A\phi)\,.
$$
\item If $\H$ is a complex Hilbert space with scalar product 
(conjugate-linear with respect to the first variable)
$\langle\cdot,\cdot\rangle$, then $C_\H:\H\to\H'$ denotes the
conjugate-linear isomorphism
defined by $$(C_\H\, \psi)(\phi):=\langle \psi,\phi\rangle\,.$$ 
\item The Hilbert
adjoint $A^*\in\L(\H_2,\H_1)$ of the densely defined linear operator 
$A\in\L(\H_1,\H_2)$ is defined as $$A^*:=
C_{\H_1}^{-1}\cdot A'\cdot C_{\H_2}\,.$$ 
\item $F$ and $*$ denote Fourier transform and convolution
respectively.
\item $H^s(\RE^n)$, $s\in\RE$, is the usual scale of
Sobolev-Hilbert spaces, i.e. $H^s(\RE^n)$ is the space of tempered
distributions with a Fourier transform which is square integrable
with respect to the measure with density $(1+|x|^2)^s$. 
\end{itemize}

\section{Extensions by a Kre\u\i n-like formula}
Given the Hilbert space
$\H$ with scalar product $\langle\cdot,\cdot\rangle$ (we denote by
$\|\cdot\|$ 
the corresponding norm and put $C\equiv C_\H$), let $A:D(A)\subseteq\H\to\H$ be a
self-adjoint operator and let $\N\subsetneq D(A)$ be a linear dense set which is 
closed with respect to the graph
norm on $D(A)$. We denote by $\HP$ the Hilbert space given by the set
$D(A)$ equipped with the scalar product $\langle\cdot,\cdot\rangle_+$
leading to the graph norm, i.e.
$$\langle\phi_1,\phi_2\rangle_{+}
:=\langle(A^2+1)^{1/2}\phi_1,(A^2+1)^{1/2}\phi_2\rangle\,.$$
We remark that in the sequel we will avoid to identify $\H_+$ with
its dual. Indeed we will use the duality map induced by the scalar
product on $\H$ (see the next section for the details).\par 
Being $\N$ closed we have $\HP=\N\oplus \N^\perp$ and we can then
consider the orthogonal projection $\pi:\HP\to \N^\perp$. 
From now on, since this gives advantages in concrete applications where
usually a variant of $\pi$ is what is known in advance, more generally
we will consider a linear map
$$
\tau:\HP\to\fh\,,\qquad\tau\in\B(\HP,\fh)\,,
$$ 
where $\fh$ is a Hilbert space with scalar product
$\langle\cdot,\cdot\rangle_\fh$ and corrsponding norm $\|\cdot\|_\fh$, such that
\begin{equation}
\text{\rm Range}\,\tau=\fh
\end{equation}
and
\begin{equation}
\overline{\text{\rm Kernel}\,\tau}=\H\,,
\end{equation}
the bar denoting here the closure in $\H$.
We put $$\N:=\text{\rm
Kernel}\,\tau\,.$$
By (2.1) one has $\fh\simeq \HP/\text{\rm Kernel$\,\tau$}
\simeq \N^\perp$ so that 
$$\HP\simeq \N\oplus\fh\,.$$
Regarding (2.2) we have the following
\begin{lemma} Hypothesis
(2.2) is equivalent to 
\begin{equation*}
\text{\rm Range}\,\tau'\cap\H'=\left\{0\right\}\,,
\end{equation*}
when one uses the embedding of $\H'$ into $\HP'\supseteq\text{\rm
Range}\,\tau'$ 
given by the map 
$\phi\mapsto \langle C^{-1}\phi,\,\cdot\,\rangle$. \end{lemma} 
\begin{proof} Defining as usual the annihilator of $\N$ by 
$$
\N^0:=\left\{\lambda\in\HP'\ :\ \forall\,\phi\in\N,\quad\lambda(\phi)=0\right\}
$$ 
one has that denseness of $\N$ is equivalent to
$$
\N^0\cap\H'=\left\{0\right\}\,.
$$
Since $\overline{\text{\rm
Range}\,\tau'}=\N^0$ the proof is concluded if the
range of $\tau'$ is closed. This follows from the closed range theorem
since the range of $\tau$ is closed by the surjectivity hypothesis.
\end{proof}
Being $\rho(A)$ the resolvent set of $A$, we define $R(z)\in \B(\H,\HP)$,
$z\in\rho(A)$, by
$$
R(z):=(-A+z)^{-1}
$$
and we then introduce, for any $z\in\rho(A)$, the two linear
operators $\breve G(z)\in
\B(\H,\fh)$ and $G(z)\in\B(\fh,\H)$ by
$$
\breve G(z):=\tau\cdot R(z)\,,\qquad 
G(z):=\breve G(\z*)^*\,.
$$
By (2.2) one has 
\begin{equation}
\text{\rm Range}\,G(z)\cap D(A)=\left\{0\right\}\,,
\end{equation}
and, as an immediate consequence of the first resolvent 
identity for $R(z)$  (see \cite{[P1]}, Lemma 2.1) 
\begin{equation}
(z-w)\,R(w)\cdot G(z)=G(w)-G(z)\,.
\end{equation}
These relations imply
\begin{equation}
\text{\rm
Range}\,(G(w)-G(z))\subseteq D(A)
\end{equation} 
and
$$
\text{\rm
Range}\,(G(w)+G(z))\cap D(A)=\left\{0\right\}\,.
$$
By \cite{[P1]} (combining Theorem 2.1,
Proposition 2.1, Lemma 2.2, Remarks 2.10, 2.12 and 2.13) one then obtains the
following 
\begin{theorem} Given $z_0\in\C\backslash\RE$ define 
$$
\G*:=\frac{1}{2}\,(G(z_0)+G(\z*_0))\,\qquad 
G_\diamond:=\frac{1}{2}\,(G(z_0)-G(\z*_0))\,
$$
and, given then any self-adjoint operator $\Theta:
D(\Theta)\subseteq\fh\to\fh$, define  
$$
R_\Theta(z):=R(z)+G(z)\cdot\left(\Theta+\Gamma(z)\right)^{-1}\cdot\breve
G(z)\,,
\qquad z\in W_\Theta\cup\C\backslash\RE\,,
$$
where
$$
\Gamma(z):=\tau
\cdot\left(\G*-G(z)\right)
$$
and
$$
W_\Theta:=\left\{\,\lambda\in\RE\cap \rho(A)\ :\ 
0\in\rho(\Theta+\Gamma(\lambda))\,\right\}\,.
$$
Then $R_\Theta$ is the resolvent of the self-adjoint extension of $A_\N$
defined by 
\begin{align*}
D(A_\Theta):=
\left\{\,\phi\in\H\, :\, \phi=
\ph*+\G*Q_\phi,\right.&\\
\left. \ph*\in D(A),\ Q_\phi\in D(\Theta),\ 
\tau\ph*=\Theta\,Q_\phi\,\right\}&\,,
\end{align*}
$$
A_\Theta\,\phi:=A\,\ph*+\text{\rm Re}(z_0)\,\G*Q_\phi+i\,\text{\rm
Im}(z_0)\,G_\diamond Q_\phi\,.
$$
\end{theorem}
\begin{proof} Here we just give the main steps of the proof refering
to \cite{[P1]}, \S 2, for the details. One starts writing the presumed
resolvent of an extension $\wtilde A$ of $A_\N$ as 
$$
\wtilde R(z)=R(z)+B(z)\cdot\tau\cdot R(z)\equiv R(z)+B(z)\cdot\breve G(z)\,, 
$$
where $B(z)\in \B(\fh,\H)$ has to be determined. 
Self-adjointness requires $\wtilde R(z)^*=\wtilde R(\z*)$ or, equivalently, 
\begin{equation}
G(\z*)\cdot B(z)^*=B(\z*)\cdot\breve G(\z*)\ .
\end{equation}
Therefore posing $B(z)=G(z)\cdot \Lambda(z)$, where
$\Lambda(z)\in\B(\fh)$, 
(2.6) is equivalent
to 
\begin{equation}
\Lambda(z)^*=\Lambda(\z*)\ .
\end{equation}
The resolvent identity 
\begin{equation}
(z-w)\,\wtilde R(w)\wtilde R(z)=\wtilde R(w)-\wtilde R(z)
\end{equation}
is then equivalent to  
\begin{equation}
\Lambda(w)-\Lambda(z)=
(z-w)\,\Lambda(w)\cdot\breve G(w)\cdot
G(z)\cdot\Lambda(z)\ .
\end{equation}
Suppose now that there 
exist a (necessarily closed) operator
$$\Gamma(z):D\subseteq\fh\to\fh 
$$
and an open set $Z\subseteq\rho(A)$, invariant with respect to complex conjugation, such that  
$$
\forall\, z\in Z,\qquad\Gamma(z)^{-1}=\Lambda(z)\ .
$$
Then (2.9) forces $\Gamma(z)$ to satisfy the relation
\begin{equation}
\Gamma(z)-\Gamma(w)=(z-w)\,\breve G(w)\cdot
G(z)\,,
\end{equation}
whereas (2.7), at least in the case $\Gamma(z)$ is densely defined, and
has a bounded inverse given by
$\Lambda(z)$ 
as we are pretending, is equivalent to
\begin{equation}
\Gamma(z)^*=\Gamma(\z*)\,.
\end{equation}  
By \cite{[P1]}, Lemma 2.2, for any self-adjoint $\Theta$, the linear
operator 
$$\Theta+\tau
\cdot\left(\G*-G(z)\right)$$ satisfies (2.10), (2.11) and, by
\cite{[P1]}, Proposition 2.1, has
a bounded inverse for any $z\in W_\Theta\cup\C\backslash\RE$ 
(at this point hypothesis (2.1) is used).
Therefore (see the proof of Theorem 2.1 in \cite{[P1]})
$$
R_\Theta(z):=R(z)+G(z)\cdot\left(\Theta+\Gamma(z)\right)^{-1}\cdot\breve
G(z)
$$
is the resolvent of a self-adjoint operator $A_\Theta$ (here hypotheses (2.2) is
needed). For any $z\in W_\Theta\cup\C\backslash\RE$ one has  
\begin{equation}
D(A_\Theta)=\left\{\,\phi\in\H\, :\, \phi=
\phi_z+G(z)\cdot(\Gamma+\Theta(z))^{-1}\cdot\tau\,\phi_z,\, \phi_z\in
D(A)\,\right\}\,,
\end{equation}
\begin{equation}
(-A_\Theta+z)\phi=(-A+z)\phi_z\,,
\end{equation}
the definition of $A_\Theta$ being $z$-independent thanks to resolvent
identity (2.8). Being $G(z)$ injective, (2.3) and (2.5) imply
$$
\phi_w+G(w)Q_1=\phi_z+G(z)Q_2\quad\Rightarrow\quad  Q_1=Q_2
$$
and so the definition 
$$
Q_\phi:=(\Theta+\Gamma(z))^{-1}\cdot\tau\,\phi_z
$$
is $z$-independent. Therefore
any $\phi\in D(A_\Theta)$ can be equivalently 
re-written as 
$$
\phi=\phi_z+G(z)Q_\phi\,,
$$
where $Q_\phi\in D(\Theta)$ and 
$$
\tau\phi_z=\Theta\,Q_\phi+\Gamma(z)Q_\phi\,.
$$
This implies, for any $\phi\in D(A_\Theta)$,
\begin{align*}
\phi=&\frac{1}{2}\,\left(\phi_{z_0}+G(z_0)Q_\phi+
\phi_{\z*_0}+G(\z*_0)Q_\phi\right)
\equiv\ph*+\G*Q_\phi\,,\\
\tau\ph*\equiv&\frac{1}{2}\,\tau(\phi_{z_0}+\phi_{\z*_0})
=\Theta\,Q_\phi+\frac{1}{2}\,(\Gamma(z_0)Q_\phi+\Gamma(\z*_0)Q_\phi)=\Theta\,Q_\phi\,,\\
A_\Theta\phi=&\frac{1}{2}\,(A\phi_{z_0}+z_0G(z_0)Q_\phi
+A\phi_{\z*_0}+\z*_0G(\z*_0)Q_\phi)\\
\equiv& 
A\ph*+\text{\rm Re}(z_0)\,\G*Q_\phi+i\,\text{\rm
Im}(z_0)\,G_\diamond Q_\phi\,.
\end{align*}
Conversely any $\phi=\ph*+\G*Q_\phi$, $\ph*\in D(A)$, $\Theta\,
Q=\tau\ph*$, admits the decomposition 
$\phi=\phi_z+G(z)\cdot (\Theta+\Gamma(z))^{-1}\cdot\tau\phi_z$, where
$$
\phi_z:=\ph*+(\G*-G(z))Q_\phi\,.
$$ 
Note that $\phi_z\in D(A)$ by (2.5) and $\tau\phi_z=(\Theta+\Gamma(z))Q_\phi$.
\end{proof}
\begin{remark} The results quoted in the previous theorem are
consequences of an alternative version of Kre\u\i n's resolvent formula. 
The original one was obtained in \cite{[K1]}, \cite{[K2]},  
\cite{[S]} for the cases where dim$\,\K_\pm=1$, dim$\,\K_\pm<+\infty$,
dim$\,\K_\pm=+\infty$ respectively; also see \cite{[DM]}, 
\cite{[GMT]}, \cite{[KK]} for more recent formulations. In standard
Kre\u\i n's formula (usually written with $z_0=i\equiv\sqrt{-1}\,$)
the main ingredient is the orthogonal projection $P:\H\to \K_+$ whereas
we used, exploiting the a priori knowledge of the self-adjoint operator $A$, the
map $\tau$, 
which plays the role of the orthogonal projection
$\pi: \HP\to \N^\perp$. Thus the knowledge of $A_\N^*$ is not needed. 
The version given in
\cite{[P1]} allows $\tau$ to be not surjective and $\fh$ can be a
Banach space; the use of
the map $\tau$ simplifies the exposition and makes easier to work out
concrete applications. Indeed, as we already said, frequently what is explicitely known 
is the map $\tau$ and $\N$ is then simply defined as its
kernel: see the many examples in \cite{[P1]} where $\tau$
is the trace (restriction) map along some null subset of
$\RE^n$ and $A$ is a (pseudo-)differential operator. Moreover this approach 
allows a natural formulation in terms of the boundary
condition $\tau\ph*=\Theta\,Q_\phi$. Note that, since $\G* Q_\phi\in D(A)$ if and
only if $Q_\phi=0$, once the reference point $z_0$ has been chosen,
the decomposition $\phi=\ph*+\G* Q_\phi$ of a generic element $\phi$
of $D(A_\Theta)$ by a regular part
$\ph*\in D(A)$ and a singular one $\G* Q_\phi\in \H\backslash D(A)$ is univocal. 
\end{remark}
\begin{remark}  As regards the definition of $R_\Theta(z)$, the one
given in the theorem above is not the only
possible definition of the operator $\Gamma(z)$. 
Any other not
necessarily bounded, densely defined operator 
satisfying 
$$
\Gamma(z)-\Gamma(w)=(z-w)\,\breve G(w)\cdot G(z)\,,
$$
$$
\Gamma(\z*)= \Gamma(z)^*
$$
and such that $\Theta+\Gamma(z)$ is boundedly invertible would
suffice; moreover hypothesis (2.1) is not necessary 
(see \cite{[P1]}, Theorem 2.1); 
note that, once $\Theta$ is given, $\Gamma(z)$ univocally defines
$(-A_\Theta+z)^{-1}$ and hence $A_\Theta$ itself. 
For alternative
choices of $\Gamma(z)$ we refer to \cite{[P1]}; also see \cite{[P2]} where it is
shown how, 
under the hypotheses Kernel$\,A=\left\{0\right\}$ and $||\tau\phi||_\fh\le
c\, \|A\phi\|$,
it is always possible to take $z_0=0$ in Theorem 2.1 (at the expense of
having then $\ph*$ in the completion of $D(A)$ 
with respect to the norm $\phi\mapsto \|A\phi\|\,$). 
However we remark that any different choice (either of $z_0$ or of the
operator $\Gamma(z)$ itself) does not change the 
family of extensions as a whole. 
\end{remark}
\begin{remark} In the case $A$ has a non-empty real resolvent set, by \cite{[P1]}, Remark 2.7, if in Theorem 2.1 one consider only the
sub-family of extensions in which the $\Theta$'s have
bounded inverses, then one can take $z_0\in\RE\cap\rho(A)$. More
generally one can take $z_0\in W_\Theta$ independently of the invertibility of $\Theta$; however this could give rise to implicit conditions 
(related to the location of the spectrum of $A_\Theta$) on the
choice of $z_0$. \par
\end{remark}
\section{Extensions by Additive Perturbations}
We define the pre-Hilbert space $\tilde\HM$ as the set $\H$
equipped with the scalar product 
$$
\langle\phi_1,\phi_2\rangle_-:=
\langle (A^2+1)^{-1/2}\phi_1,(A^2+1)^{-1/2}\phi_2\rangle\,.
$$
We denote then by $\HM$ the Hilbert space given by the completion of
$\tilde\HM$.
We will avoid to identify $\HP$ and $\HM$ with their duals; indeed,
see Lemma 3.1 below, we will identify $\HP'$ with $\HM$.\par
As usual
$\H$ will be treated as a (dense) subspace of $\HM$ by means of the
canonical embedding $$I_-:\H\to\HM$$ which associates to $\phi$ the set
of all the Cauchy sequences converging to $\phi$. 
Considering also the canonical embedding (with dense range) 
$$
I_+:\H'\to\HP'\,,\qquad I_+\lambda(\phi):=\langle
C^{-1}\lambda,\phi\rangle\,,
$$
we can then define the conjugate linear operator 
$$
C_-:\HP'\to\HM
$$
as the unique bounded extension of 
$$
I_-\cdot C^{-1}\cdot I_+^{-1}:I_+(\H')\subseteq\HP'\to\HM\,.
$$
Analogously we define the conjugate linear operator
$$
C_+:\HM\to\HP'
$$
as the unique bounded extension of 
$$
I_+\cdot C\cdot I_-^{-1}:I_-(\H)\subseteq\HM\to\HP'\,.
$$
These definitions immediately lead to the following
\begin{lemma} One has
$$C_+=C_-^{-1}\,,\qquad C_-=C_+^{-1}\,,$$
so that
$$
\HP'\,\simeq\,\HM\,.
$$
\end{lemma}
We will denote by 
$$
(\cdot,\cdot):\HM\times\HP\to\C\,,\qquad (\vp,\phi):=C_+\vp(\phi)
$$
the pairing between $\HM$ and $\HP$. It is nothing else that the
extension of the scalar product of $\H$, being  
$$
(I_-\phi_1,\phi_2)=\langle\phi_1,\phi_2\rangle\,.
$$
We consider now the linear operator 
$$I_-\cdot A:\HP\subseteq\H\to\HM\,.$$
Since
$$
\|(A^2+1)^{-1/2}A\phi\|\le\|\phi\|\,,
$$
the operator $I_-\cdot A$ has an unique extension 
$$\A:\H\to\HM\,,\qquad\A\in\B(\H,\HM)\,.$$
\begin{lemma} Let $A':\H'\to\HP'$ be the adjoint of the linear operator
$A$ when viewed as an element of $\B(\HP,\H)$. Then one has 
$$
\A=C_-\cdot A'\cdot C\,.
$$
\end{lemma}
\begin{proof} Being $I_-$ injective, by continuity and density the
thesis follows from the identity
$$
A=A^*\equiv C^{-1}\cdot A'\cdot C\,.
$$
\end{proof}
\begin{remark} If we use the symbol $A_+$ to denote the linear
operator $A$ when we consider it as an element of $\B(\HP,\H)$, and if
we use $C_-$ as a substitute of $C^{-1}_{\HP}$, then by Lemma 3.2 and
a slight abuse of notations we can write
$$
\A=A_+^*\,.
$$
By the same abuse of notations we define $\tau^*\in\B(\fh,\HM)$ by
$$
\tau^*:=C_-\cdot\tau'\cdot C_\fh\,.
$$
\end{remark} 
Now we can reformulate Theorem 2.1 in terms of additive perturbations:
\begin{theorem}  
Define
$$
D(\wtilde A):=\left\{\,\phi\in\H\,:\,\phi=\ph*+\G* Q_\phi,\ \ph*\in
D(A),\ Q_\phi\in\fh\,\right\}\,,
$$
$$
\wtilde A:D(\wtilde A)\to\HM\,,\qquad \wtilde A:=\A+T\,,
$$
where
$$
T:D(\wtilde A)\to\HM\,, \qquad T\phi:=\tau^* Q_\phi\,. 
$$
Then the linear operator $\wtilde A$ is $\H$-valued and coincides with 
$A_\Theta$ when restricted to $D(A_\Theta)$, i.e. when a boundary
condition of the kind $\tau\ph*=\Theta\,Q_\phi$ holds for some self-adjoint operator
$\Theta$. Therefore, posing 
$T_\Theta:={T\,}_{\left|D(A_\Theta)\right.}$, one has 
$$
A_\Theta:D(A_\Theta)\to\H\,,\qquad A_\Theta:=\A+
T_\Theta\,,$$
and, in the case $\Theta$
has a bounded inverse, 
$$
A_\Theta:D(A_\Theta)\to\H\,,\qquad A_\Theta\phi=\A\phi+V_\Theta\ph*\,,$$
where
$$
V_\Theta:\HP\to\HM\,,\qquad (V_\Theta\phi_1,\phi_2)
:=\langle\Theta^{-1}\tau\phi_1,\tau\phi_2\rangle_\fh\,.
$$
\end{theorem}
\begin{proof} By the definition of $\A$, $\tau^*$ and $\G*$ one has, for any
$\phi\in D(\wtilde A)$,
\begin{align*}
\A\phi=&I_-\cdot A\ph*+C_-\cdot A'\cdot C\cdot \G* Q_\phi\\
=&I_-\cdot A\ph*+\frac{1}{2}\,C_-\cdot A'\cdot
R(\z*_0)'\cdot\tau'\cdot C_\fh\, Q_\phi\\
&+\frac{1}{2}\,C_-\cdot A'\cdot
R(z_0)'\cdot\tau'\cdot C_\fh\, Q_\phi\\
=&I_-\cdot (A\ph*+\text{\rm Re}(z_0)\,\G*Q_\phi+i\,\text{\rm
Im}(z_0)\,G_\diamond Q_\phi)-T Q_\phi \,.
\end{align*}
The proof is then concluded by Theorem 2.1.
\end{proof}
\begin{remark} In the case $0\in\rho(A)$ and $\Theta$ is boundedly
invertible, by Theorem 2.1 and Remark 2.4 (taking $z_0=0$) one can
define $A_\Theta$ either by $A_\Theta\phi:=A\ph*$ or, equivalently, by 
$$
A_\Theta^{-1}=A^{-1}+ G\cdot\Theta^{-1}\cdot \breve G\,,
$$
where $G:=G(0)$, $\breve G:=\breve G(0)$.
Since, for any 
$\phi_1,\phi_2\in\H$, one has 
\begin{align*}
&\langle \A^{-1}\cdot V_\Theta\cdot
A^{-1}\phi_1,\phi_2\rangle =(V_\Theta A^{-1}\phi_1, A^{-1}\phi_2)\\
=&
\langle\Theta^{-1}\tau\cdot A^{-1}\phi_1,\tau\cdot
A^{-1}\phi_2\rangle_\fh=
\langle\Theta^{-1}\breve G\phi_1,\breve G\phi_2\rangle_\fh\\
=&
\langle G\cdot\Theta^{-1}\cdot\breve G\phi_1,\phi_2\rangle\,,
\end{align*}
the self-adjoint extension $A_\Theta$ could be defined directly in
terms of $V_\Theta$ by
$$
A_\Theta^{-1}=A^{-1}+\A^{-1}\cdot V_\Theta\cdot A^{-1}\,.
$$
This reproduces the formulae appearing in \cite{[AKK]},
Lemma 2.3, where however no additive representaion of the extension $A_\Theta$
is given, and in \cite{[K]} where an additive representaion is
obtained only when $\N$ is closed in $D(A^{1/2})$.
\end{remark}
\section{The connection with von Neumann's Theory}
In this section we explore the connection between the
results given in the previous sections and von Neumann's theory of
self-adjoint extensions \cite{[N]}. 
Such a theory (see e.g. \cite{[F]}, \S 13, for a very compact exposition)
tells us that
$$
D(A^*_\N)= \N\oplus
\K_+\oplus \K_-\,,\quad
A_\N^*(\phi_0+\phi_++\phi_-)=A\phi_0+i\phi_+-i\phi_-\,,
$$ 
the direct sum decomposition being orthogonal with
respect to the graph inner product of $A^*_\N$; any self-adjoint
extension $A_U$ of $A_\N$ is then obtained
by restricting $A_\N^*$ to a subspace of the kind $\N\oplus \text{\rm
Graph}\,U$, where
$U:\K_+\to K_-$ is unitary. \par For simplicity in the next theorem we will consider only the case
$z_0=i$ and we put $G_\pm:=G(\pm i)$ and $\Gamma:=\Gamma(i)$.
\begin{theorem} Let $\wtilde A=\A+T$ as defined in Theorem 3.4. Then 
$$\wtilde
A=A_\N^*\,.
$$
The linear operator 
$$G_\pm:\fh\to
\K_\pm$$ is a continuos bijection which becomes unitary when one puts on $\fh$ the
scalar product
$$
\langle Q_1,Q_2\rangle_\Gamma:=
\langle\sqrt{-i\Gamma}\,Q_1,\sqrt{-i\Gamma}\,Q_2\rangle_\fh\,.
$$
The linear operator 
$$U:\K_+\to \K_-\,,\qquad
U:=-\,G_-\cdot(\uno+2(\Theta-\Gamma)^{-1}\cdot\Gamma)\cdot G_+^{-1}
$$
is unitary and the corresponding von Neumann's extension $A_U$
coincides with the self-adjoint operator $A_\Theta$ defined in
Theorems 2.1 and 3.4.
\end{theorem}
\begin{proof} By the definition of $\breve G_\pm\equiv\breve G(\pm i)$ one has 
$$
\text{\rm Range}\, (-A_\N\pm i)=\text{\rm Kernel}\,\breve G_\pm
$$
and so, since
$$
\K_\pm=\text{\rm Range}\,(-A_\N\mp i)^\perp
$$
and
$$
{\rm Range}\,{G_\pm}^\perp=\text{\rm Kernel}\,\breve G_\mp\,,
$$
in conclusion there follows
$$
{\rm Range}\,G_\pm=\K_\pm
$$
if and only if Range$\,G_\pm$ is closed. By the closed range theorem
Range$\,G_\pm$ is 
closed if and only if 
the Range$\,\breve G_\pm$ is closed, and this is equivalent to the
range of $\tau$ being closed. 
Being $\tau$ surjective, $G_\pm$ is injective with a closed range and
so 
$$G_\pm:\fh\to
\K_\pm$$ is a bijection.\par
By von Neumann's theory we know that any $\phi\in D(A_\N^*)$
can be univocally decomposed as 
$$
\phi=\phi_0+\phi_++\phi_-\,,\qquad\phi_0\in \N,\ \phi_\pm\in \K_\pm\,,
$$
i.e.
$$
\phi=\phi_0+G_+Q_++G_-Q_-\,,\qquad\phi_0\in \N,\ Q_\pm\in \fh\,.
$$
The above decomposition can be then rearranged as
\begin{align*}
\phi=&\phi_0+\frac{1}{2}\,(G_+-G_-)Q_++\frac{1}{2}\,(G_++G_-)Q_+\\
&\quad
+\frac{1}{2}\,(G_--G_+)Q_-+\frac{1}{2}\,(G_-+G_+)Q_-\\
=&\phi_0+\frac{1}{2}\,(G_--G_+)(Q_--Q_+)+\G*(Q_-+Q_+)\,.
\end{align*}
By (2.4) one has
\begin{equation}
G_\mp-G_\pm=\pm 2i\,R(\mp i)\cdot G_\pm\,.
\end{equation}
Since the scalar product of $\HP$ can be equivalently  written
as 
$$
\langle\phi_1,\phi_2\rangle_+=\langle(-A+i)\phi_1,(-A+i)\phi_2\rangle\,,
$$
one has
$$
G_--G_+=2i\,R(-i)\cdot R(-i)^*\cdot\tau^*=2i\,\tau^*\,.
$$
This implies, since Range$\,G_+$ is closed,
$$
\text{Range}\,(G_--G_+)=\text{Range}\,\tau^*=\text{Kernel}\,{\tau\,}^\perp\,.
$$
Thus, being $\HP=\N\oplus \N^\perp$, the vector
$$
\phi_0+\frac{1}{2}\,(G_--G_+)(Q_--Q_+)
$$
is a generic element of $D(A)$ and we have shown that
$D(\wtilde A)=D(A^*_\N)$. It is then straighforward to check that $\wtilde
A=A_\N^*$.\par
By (4.1) one has
$$
\Gamma=\pm\frac{1}{2}\,\tau\cdot(G_\mp-G_\pm)=i\,\tau\cdot
R(\mp i)\cdot G_\pm=i\,\breve G_\mp\cdot G_\pm=i\,G_\pm^*\cdot G_\pm\,.
$$
This implies 
$$
\|G_\pm Q\|=\|\sqrt{-i\Gamma}\,Q\|_\fh\,,
$$
thus $U=-\,G_-\cdot(\uno+2(\Theta-\Gamma)^{-1}\cdot\Gamma)\cdot
G_+^{-1}$ is isometric if and only if 
$$
\forall\, Q\in\fh\,,\qquad
\|\sqrt{-i\Gamma}\cdot\wtilde U\,Q\|_\fh=\|\sqrt{-i\Gamma}\,Q\|_\fh\,,
$$
where $\wtilde U:=G_-^{-1}\cdot U\cdot G_+$. By using the identities $\Gamma^*=-\Gamma$
and 
\begin{equation}
(\Theta-\Gamma)^{-1}-(\Theta+\Gamma)^{-1}
=2\,(\Theta+\Gamma)^{-1}\cdot\Gamma\cdot(\Theta-\Gamma)^{-1}\,,
\end{equation}
one has
\begin{align*}
&i\,\|\sqrt{-i\Gamma}\cdot(\uno+2(\Theta-\Gamma)^{-1}\cdot \Gamma)\,Q\|_\fh\\
=&\langle\Gamma\, Q+2\Gamma\cdot(\Theta-\Gamma)^{-1}\cdot \Gamma\,Q,
Q+2(\Theta-\Gamma)^{-1}\cdot \Gamma\,Q\rangle
\\=&\langle\Gamma Q,Q\rangle
+2\langle\Gamma\,Q,(\Theta-\Gamma)^{-1}\cdot \Gamma\,Q\rangle
+2\langle\Gamma\cdot(\Theta-\Gamma)^{-1}\cdot \Gamma\,Q,Q\rangle\\
&+4\langle\Gamma\cdot(\Theta-\Gamma)^{-1}\cdot \Gamma\,Q,
(\Theta-\Gamma)^{-1}\cdot \Gamma\,Q\rangle\\
=&\langle\Gamma Q,Q\rangle
+2\langle\Gamma\,Q,((\Theta-\Gamma)^{-1}-(\Theta+\Gamma)^{-1})\cdot \Gamma\,Q\rangle\\
&-4\langle \Gamma\,Q,
(\Theta+\Gamma)^{-1}\cdot\Gamma\cdot(\Theta-\Gamma)^{-1}\cdot \Gamma\,Q\rangle\\
=&\langle\Gamma Q,Q\rangle=i\,\|\sqrt{-i\Gamma}\,Q\|_\fh\,,
\end{align*}
and so $U$ is an isometry. By again using identity (4.2) one can
check that $U$ has an inverse defined by
$$
U^{-1}:=-\,G_+\cdot(\uno-2(\Theta+\Gamma)^{-1}\cdot\Gamma)\cdot G_-^{-1}\,.
$$
Thus $U$ is unitary. Let us now take $G_-Q_-=U G_+Q_+$. Then
$$
-2(\Theta-\Gamma)^{-1}\cdot\Gamma Q_+=Q_-+Q_+
$$
and so $Q_-+Q_+\in D(\Theta)$ and
$$
\tau\left(\phi_0+\frac{1}{2}\,(G_--G_+)(Q_--Q_+)\right)
\equiv\Gamma(Q_--Q_+)=\Theta(Q_-+Q_+)\,.
$$  
\end{proof}
\begin{remark} Note that when $\Theta$ is bounded, in the
previuos theorem one can re-write the unitary $U$  as
$$
U=-G_-\cdot(\Theta-\Gamma)^{-1}\cdot(\Theta+\Gamma)\cdot G_+^{-1}\,.
$$
Being $\Theta$ always bounded when dim$\K_\pm=n$, the previous 
theorem gives an analogue of Theorem 3.1.2 in \cite{[AK]} avoiding
however the use of an admissible matrix $R$ (see \cite{[AK]}, definition 3.1.2).  
\end{remark}
The previous theorem has the following converse:
\begin{theorem} Let $A_U$ be a self-adjoint
extension of $A_\N$ as given by von Neumann's theory. Suppose that 
$D(A_U)\cap D(A)=\N$ and let 
$U_A:=(-A+i)\cdot(-A-i)^{-1}$ be the Cayley
transform of $A$. Then the set  
$$
D(\Theta):=\text{\rm Range}\ G_-^{-1}\cdot(U+U_A)
$$ 
is dense,  $$
\Theta:D(\Theta)\subseteq\fh\to\fh\,,\qquad
\Theta:=i\, \breve G_+\cdot(U-U_A)\cdot(U
+U_A)^{-1}\cdot G_-\,,
$$ 
is self-adjoint and the corresponding self-adjoint operator
$A_\Theta$, defined in Theorems 2.1 and 3.4, coincides with $A_U$.   
\end{theorem}
\begin{proof} By (4.1) one has
$$
G_-\cdot G_+^{-1}=\uno + 2i\,R(-i)=U_A\,.
$$
Thus, by inverting the relation 
$U=-\,G_-\cdot(\uno+2(\Theta-\Gamma)^{-1}\cdot\Gamma)\cdot G_+^{-1}$
given in the previous theorem, one obtains
\begin{align*}
\Theta=&\Gamma\cdot(G_-^{-1}\cdot U\cdot
G_+-\uno)\cdot(G_-^{-1}\cdot U\cdot G_+
+\uno)^{-1}\\
=&\Gamma\cdot G_-^{-1}\cdot(U
-G_-\cdot G_+^{-1})\cdot(U
+G_-\cdot G_+^{-1})^{-1}\cdot G_-\\
=&\Gamma\cdot G_-^{-1}\cdot(U-U_A)\cdot(U
+U_A)^{-1}\cdot G_-\,.
\end{align*}
Since $U=-U_{A_U}$ and
$1\notin\sigma_p(U_A\cdot U_{A_U}^{-1})$ 
if and only if $D(A_U)\cap D(A)=\N$ (see
e.g. \cite{[GMT]}, Lemma 1), the range of $U+U_A$ is dense and thus
$\Theta$ is densely defined as $G_-$ is a continuos bijection. By
(4.1) one has
$$
\Gamma\cdot G_-^{-1}=i\tau\cdot R(i)\equiv i\breve G_+
$$  
and so, since $\breve G_+^*=G_-$ and $G^*_-=\breve G_+$, $\Theta$ is
self-adjoint if and only if
$$
(U^*+U_A^*)\cdot(U-U_A)=-(U^*
-U_A^*)\cdot(U
+U_A)\,.
$$ 
Such an equality is then an immediate conseguence of the unitarity of both
$U$ and $U_A$.
\end{proof}
\begin{corollary}
$\wtilde A=\A+T$ as defined in Theorem 3.4 coincides with a 
self-adjoint extension $\widehat A$ of $A_\N$ such that $D(\widehat A)\cap
D(A)=\N$ if and only if the boundary
condition $\tau\ph*=\Theta\,Q_\phi$ holds for some self-adjoint operator
$\Theta:D(\Theta)\subseteq\fh\to\fh$. 
\end{corollary}
\section{Examples}
\begin{example} {\it Finite rank perturbations.}  Suppose
dim$\,\K_\pm=n$, so that $\fh\simeq\C^n$ and
$\tau\in\B(\HP,\C^n)$. Then necessarily 
$$
\tau:\HP\to\C^n\,,\qquad\tau\phi=\left\{(\vp_j,\phi)\right\}_1^n\,,
$$
with $\vp_1,\dots,\vp_n\in\HM$. 
Hypotheses (2.1) and
(2.2) correspond to 
$$
\exists\,\phi_1,\dots,\phi_n\in\HP\quad \text{\rm s.t.}\quad
(\vp_i,\phi_j)=\delta_{ij}\,,
$$
and
$$
\sum_{j=1}^{n}c_j\,\vp_j\in\H\quad\text{\rm iff}\quad c_1=\dots =c_n=0\,.
$$ 
Considering then an Hermitean invertible matrix $\Theta=(\theta_{ij})$ with inverse 
$\Theta^{-1}=(t_{ij})$, by Theorem 3.4 one can define the self-adjoint operator  
$$A_\Theta\phi:=\A\phi+\sum_{i,j=1}^{n}t_{ij}(\vp_i,\ph*)\vp_j$$
with 
\begin{align*}
D(A_\Theta):=\left\{\,\phi\in\H\,:\,\phi=\ph*+\sum_{j=1}^{n} Q_j R_\star
\vp_j,\right.&\\
\left. \ph*\in
D(A),\ Q\in\C^n,\ (\vp_i,\ph*)=\sum_{j=1}^{n}\theta_{ij}Q_j\,\right\}&\,,
\end{align*}
where
$$
R_\star:=\frac{1}{2}\,(\hat R(z_0)+\hat R(\z*_0))\,,$$
$$
\hat R(z):\HM\to\H\,,\qquad\langle \hat R(z)\vp,\phi\rangle
:=(\vp,R(\z*)\phi)
\,.$$
According to Theorem 2.1 its resolvent is given by
$$
(-A_\Theta+z)^{-1}=(-A+z)^{-1}+\sum_{i,j=1}^{n}(\Theta+\Gamma(z))^{-1}_{ij}
\hat R(z)\vp_i\, \hat R(\z*)\vp_j\,,
$$
where
$$
\Gamma(z)_{ij}=\frac{1}{2}\,(\vp_i,(\hat R(z_0)+\hat R(z_0)-2\hat R(z))\vp_j)\,.
$$
The operator $A_\Theta$ above coincides with a generic finite rank
perturbation of the self-adjoint operator $A$ as defined in
\cite{[AK]}, 
\S 3.1. In order to realize that the resolvent written above (in the case
$z_0=i$) is the same given there, the identity 
$$
\frac{1}{2}\, \left(R(i)+R(-i)-2R(z)\right)=(1+zA)\cdot(A-z)^{-1}\cdot(A^2+1)^{-1}
$$
has to be used.\par
The previous construction can be applied to the case of so-called
point interactions in three dimensions (see \cite{[AGHH]} and
references therein). Since in example 5.2 below we will consider the case
of infinitely many point interactions, here we just
treat the simplest situation in which only one point interection
(placed at the origin) is present. In this case we take
$A=\Delta$, $\H=L^2(\RE^3)$, $\HP=H^2(\RE^3)$, $\HM=H^{-2}(\RE^3)$,
and $\vp=\delta_0$. Therefore $\tau$ is simply the
evaluation map at the origin 
$$\tau:H^2(\RE^3)\to\C\,,\qquad\tau\phi=\phi(0)\,,$$ 
and we have the
family of self-adjoint operators $\Delta_\theta$, 
$\theta\in\RE\backslash\left\{0\right\}$, defined as (we take $z_0=i$)
$$
\Delta_\theta\phi:=\Delta\phi+\theta^{-1}\ph*(0)\,\delta_0
$$
on the domain
\begin{align*}
D(\Delta_\theta):=\left\{\phi\in L^2(\RE^3)\, :\,
\phi=\ph*+Q\MG_\star,\right.&\\
\left.\ph*\in H^2(\RE^3),\, Q\in\C,\, \ph*(0)=\theta\, Q\right\}&\,,
\end{align*}
where
$$
\MG_\star(x)=\cos \frac{|x|}{\sqrt 2}\quad \frac{e^{-|x|/\sqrt 2}}{4\pi|x|}\,.
$$
This reproduces the family given in \cite{[AK]}, \S 1.5.1, and coincides with the 
family $\Delta_\alpha$ given in \cite{[AGHH]},
\S I.1.1, when one takes $\alpha=\theta-({4\pi\sqrt 2}\,)^{-1}$. The case 
$\alpha=-(4\pi\sqrt 2\,)^{-1}$ can be then recovered by directly using Theorem 3.4 in the
case $\theta=0$.
\end{example}
\begin{example} {\it Infinite rank perturbations.} Suppose
dim$\,\K_\pm=+\infty$. Then (we suppose $\H$ is separable)
$\fh\simeq\ell^2(\NA)$, $\tau\in\B(\HP,\ell^2(\NA))$ and necessarily 
$$
\tau:\HP\to\ell^2(\NA)\,,\qquad\tau\phi=\left\{(\vp_j,\phi)\right\}_1^\infty\,,
$$
with $\{\vp_j\}_1^\infty\subset\HM$. The generalization of the finite rank case to
this situation is then evident. As concrete example one can consider
infinitely many point interactions in three dimensions by taking 
$A=\Delta$, $\H=L^2(\RE^3)$, $\HP=H^2(\RE^3)$, $\HM=H^{-2}(\RE^3)$ as before
and an infinite and countable set
$Y\subset\RE^3$ such that
$$
\inf_{y\not=\tilde y}\, |y-\tilde y|=d>0\,.
$$
Defining then $\vp_y:=\delta_{y}$, by \cite{[AGHH]} (see page 172) 
one has $$\tau\in\B(H^2(\RE^3),\ell^2(Y))\,,$$ 
where
$$
\tau:H^2(\RE^3)\to\ell^2(Y)\,,\qquad
\tau\phi=\left\{\phi(y)\right\}_{y\in Y}\,,
$$ and hypotheses (2.1) and
(2.2) are an immediate conseguence of the discreteness of $Y$ (see
\cite{[P1]}, example 3.4). 
By Theoren 3.4, given any invertible infinite Hermitean matrix
$\Theta=(\theta_{y\tilde y})$ with a bounded inverse 
$\Theta^{-1}=(t_{y\tilde y})$, one can then define the family of self-adjoint operators
$$
\Delta_\Theta\phi
:=\Delta\phi+\sum_{y,\tilde y\in Y}t_{y\tilde y}\,\ph*(y)\,\delta_{\tilde y}
$$
on the domain
\begin{align*}
D(\Delta_\Theta):=\left\{\phi\in L^2(\RE^3)\, :\, \phi=
\ph*+\sum_{y\in Y }Q_y\MG_\star^y,\right.&\\ 
\left.\ph*\in H^2(\RE^3),\, Q\in D(\Theta),\,
\ph*(y)=\sum_{\tilde y\in Y}\theta_{y\tilde y}Q_{\tilde y}
\right\}&\,,
\end{align*}
where $\MG_\star^{y}(x):=\MG_\star(x-y)$.
When $$\theta_{yy}=\alpha+\frac{1}{4\pi\sqrt
2}\,,\qquad \theta_{y\tilde y}=-\MG_\star(y-\tilde y)\,,\quad
y\not=\tilde y\,,$$ the self-adjoint extension 
$\Delta_\Theta$ coincides with the operator 
$\Delta_{\alpha,Y}$ given in \cite{[AGHH]}, \S III.1.1 (also see
\cite{[P1]}, example 3.4). \par
In more general situations where the set $Y$ is not discrete 
the use of the unitary isomorphism
$\fh\simeq \ell^2(\NA)$ given no advantages and, how the following
example shows, it is better to work with $\fh$
itself. \par
Let $A=\Psi$, $\H=L^2(\RE^n)$, $\HP=H^s(\RE^n)$,
$\HM=H^{-s}(\RE^n)$, where the self-adjoint pseudo-differential
operator $\Psi$ is defined by 
$$\Psi:H^s(\RE^n)\to L^2(\RE^n)\,,\qquad\Psi\phi:=F^{-1}(\psi\,F \phi)\ ,
$$
with $\psi$ is a real-valued function such that $$\frac{1}{c}\,(1+|x|^2)^{s/2}\le
1+|\psi(x)|\le c\,(1+ |x|^2)^{s/2}\,,\quad c>0\,.$$
We want now to define the self-adjoint extensions of the restriction
of $\Psi$ to functions vanishing on a $d$-set, with $0<n-d<2s$. 
A Borel set $M\subset\RE^n$ is called a $d$-set, $d\in(0,n]$, if 
$$
\exists\, c_1,\,c_2>0\ :\ \forall\, x\in M,\ \forall\,r\in(0,1),\qquad
c_1r^d\le\mu_d(B_r(x)\cap M)\le c_2r^d\ ,
$$
where $\mu_d$ is the $d$-dimensional Hausdorff measure and $B_r(x)$ is
the closed $n$-dimensional ball of radius
$r$ centered at 
the point $x$ (see \cite{[JW]}, \S 1.1, chap. VIII). 
Examples of $d$-sets are $d$-dimensional
Lipschitz submanifolds  and 
(when $d$ is not an integer) self-similar fractals of Hausdorff
dimension $d$ (see \cite{[JW]}, chap. II, example 2). We take as the
linear operator $\tau$ 
the unique continuous surjective (thus (2.1) holds true) map  
$$
\tau_M:H^s(\RE^n)\to B^{2,2}_\alpha(M)\,,\qquad\alpha=s-\,\frac{n-d}{2}
$$
such that, for $\mu_d$-a.e. $x\in M$,
$$
\tau_M\phi(x)\equiv \left\{\phi_M^{(j)}(x)\right\}_{|j|< \alpha}
=\left\{\lim_{r\downarrow 0}\,\frac{1}{\lambda_n(r)}
\int_{B_r(x)}dy\,D^j\phi(y)\,\right\}_{|j|<
\alpha}
\,,
$$
where $j\in{\mathbb Z}^n_+$, $|j|:=j_1+\dots +j_n$,
$D^j:=\partial_{j_1}\cdots\partial_
{j_n}$ and $\lambda_n(r)$ denotes the
$n$-dimensional Lebesgue measure
of $B_r(x)$. We refer to \cite{[JW]}, Theorems 1 and 3,
chap. VII, for the existence of the map $\tau_M$; obviously it coincides with the usual
evaluation along $M$ when restricted to smooth functions. 
The definition of the Besov-like space $B^{2,2}_\alpha(M)$ is quite
involved and we will 
not reproduce it here (see \cite{[JW]}, \S 2.1, chap. V). However, in the case
$0<\alpha<1$ (i.e. $2(s-1)<n-d<2s$), $B^{2,2}_\alpha(M)$ can
be alternatively defined (see \cite{[JW]}, \S 1.1, chap. V) 
as the Hilbert space of $f\in L^2(F;\mu_M)$ having finite norm
$$\|f\|^2_{B^{2,2}_\alpha(M)}:=
\|f\|^2_{L^2(M)}+\int_{|x-y|<1}d\mu_M(x)\,d\mu_M(y)\,\,
\frac{|f(x)-f(y)|^2}{|x-y|^{d+2\alpha}}\ ,
$$
where $\mu_M$ denotes the restriction of the $d$-dimensional
Hausdorff measure $\mu_d$ to the set $M$. \par
The adjoint map $\tau^*_M$ gives rive, for any $Q\in
B^{2,2}_\alpha(M)$, to the signed measure
$\nu_M(Q)\in H^{-s}(\RE^n)$ defined by
$$
(\nu_M(Q),\phi)=\langle Q,\tau_M\phi\rangle_{B^{2,2}_\alpha(M)}\,.
$$
Since $\nu_M(Q)$ has support given by the closure of $M$, hypothesis
(2.2) is always verified when the closure of $M$ has zero Lebesgue
measure. Defining then
$$
\MG^\psi_\star:=\text{\rm Re}\,F^{-1}\frac{1}{-\psi+z_0}\,, 
$$
one has
$$
G_\star: B^{2,2}_\alpha(M)\to L^2(\RE^n)\,,\qquad
G_\star Q:=\MG^\psi_\star *\nu_M(Q) \,.
$$ 
Therefore, given any self-adjoint $\Theta:D(\Theta)\subseteq
B^{2,2}_\alpha(M)\to B^{2,2}_\alpha(M)$, one has the family of self-adjoint extensions
\begin{align*}
D(\Psi_\Theta):=\left\{\phi\in L^2(\RE^n)\,:\, \phi
=\ph*+\MG^\psi_\star * \nu_M(Q_\phi)\right.&\\
\left. \ph*\in H^s(\RE^n),\, Q_\phi\in D(\Theta),\, 
\tau_M\ph*=\Theta\,Q_\phi\right\}&\,,
\end{align*}
$$
\Psi_\Theta\phi:=F^{-1}(\psi F\phi)+\nu_M(Q_\phi)
$$
(see \cite{[P1]}, example 3.6, \cite{[P2]}, \S 4, for
alternative definitions).\par
When $M$ is a compact Riemannian manifold, $\Delta_{LB}$ the
Laplace-Beltrami operator, one has
$$B^{2,2}_\alpha(M)
\simeq H^\alpha(M)=\left\{Q\in L^2(M)\,:\, (-\Delta_{LB})^{\alpha/2} Q\in L^2(M)\right\}$$ 
and 
$$
\nu_M(Q)=((-\Delta_{LB})^{\alpha}Q)\,\delta_M\,,
$$
where, for any $\wtilde Q\in H^{-\alpha}(M)\equiv H^{\alpha}(M)'$,
$$
\wtilde Q\,\delta_M(\phi):=\int_Mdv\,(-\Delta_{LB})^{-\alpha/2}
\wtilde Q\ (-\Delta_{LB})^{\alpha/2}\tau_M\phi\,,
$$
$dv$ denoting the volume element of $M$. Therefore in this case, when
$\alpha\ge 1$ (i.e. $0<n-d\le 2$), taking $\psi(k)=|k|^2$, 
$\Theta=(-\Delta_{LB})^{\alpha-1}$, one can define the
self-adjoint extension
$$
-\Delta_M\phi:=-\Delta\phi-\Delta_{LB}\cdot\tau_M\ph*\,\delta_M\,,
$$
and so the construction given here generalizes the examples given in 
\cite{[KKO]} and \cite{[Ko]}. Also see \cite{[P2]}, example 14, for an
alternative definition. 
\end{example}


\begin{thebibliography}{99}
\bibitem{[AGHH]} S. Albeverio, F. Gesztesy, R. H\o egh-Krohn, H. Holden: 
{\it Solvable Models 
in Quantum Mechanics}. Berlin, Heidelberg, New York: Springer-Verlag 1988
\bibitem{[AKK]} S. Albeverio, W. Karwowski, V. Koshmanenko: Square
Powers of Singularly Perturbed Operators. {\it Math. Nachr.} {\bf 173}
(1995), 5-24
\bibitem{[AK]} S. Albeverio, P. Kurasov: {\it Singular Perturbations of
Differential Operators}. Cambridge: Cambridge Univ. Press 2000 
\bibitem{[DM]} V.A. Derkach, M.M. Malamud: Generalized Resolvents and
the Boundary Value Problem for Hermitian Operators with Gaps. 
{\it J. Funct. Anal.} {\bf 95} (1991), 1-95 
\bibitem{[F]} W.G. Faris: {\it Self-Adjoint Operators.} Lecture Notes
in Mathematics 433. Berlin, Heidelberg, New York:
Springer-Verlag 1975
\bibitem{[GMT]} F. Gesztesy, K.A. Makarov, E. Tsekanovskii: An Addendum
to Krein's Formula. {\it J. Math. Anal. Appl.} {\bf 222} (1998),
594-606
\bibitem{[JW]} A. Jonsson, H. Wallin: Function Spaces on Subsets of $\RE^n$.
{\it Math. Reports} {\bf 2} (1984), 1-221
\bibitem{[KKO]} W. Karwowski, V. Koshmanenko, S. \^ Ota:
Schr\"odinger Operators Perturbed by Operators Related to Null
Sets. {\it Positivity} {\bf 2} (1998), 77-99
\bibitem{[Ko]} S. Kondej: Singular Perturbations of Laplace Operator
in Terms of Boundary Conditions. To appear in {\it Positivity}
\bibitem{[K]} V. Koshmanenko: Singular Operators as a Parameter of
Self-Adjoint Extensions. {\it Oper. Theory Adv. Appl.} {\bf 118}
(2000), 205-223
\bibitem{[K1]} M.G. Kre\u\i n: On Hermitian Operators with Deficiency Indices One. 
{\it Dokl. Akad. Nauk SSSR } {\bf 43} (1944), 339-342 [In Russian]
\bibitem{[K2]} M.G. Kre\u\i n: Resolvents of Hermitian Operators with Defect Index 
$(m,m)$. 
{\it Dokl. Akad. Nauk SSSR } {\bf 52} (1946), 657-660 [In Russian]
\bibitem{[KY]} M.G. Kre\u\i n, V.A. Yavryan: Spectral Shift Functions that
arise in Perturbations of a Positive Operator. {\it J. Operator
Theory} {\bf 81} (1981), 155-181
\bibitem{[KK]} P. Kurasov, S.T. Kuroda: Krein's Formula and
Perturbation Theory. Preprint, Stockholm University, 2000
\bibitem{[N]} J. von Neumann: Allgemeine Eigenwerttheorie Hermitscher 
Funktionaloperatoren. 
{\it Math. Ann.} {\bf 102} (1929-30), 49-131
\bibitem{[P1]} A. Posilicano: A Kre\u\i n-like Formula for Singular
Perturbations of Self-Adjoint Operators and Applications. 
{\it J. Funct. Anal.} {\bf 183} (2001), 109-147
\bibitem{[P2]} A. Posilicano: Boundary Conditions for Singular
Perturbations of Self-Adjoint Operators. {\it
Oper. Theory Adv. Appl.} {\bf 132} (2002), 333-346
\bibitem{[S]} Sh.N. Saakjan: On the Theory of Resolvents of a Symmetric Operator with 
Infinite Deficiency Indices. {\it Dokl. Akad. Nauk Arm. SSR} {\bf 44} (1965), 
193-198 [In Russian]

\end{thebibliography}
\end{document}